\documentclass[12pt]{article}
\usepackage{graphicx}

\usepackage[english]{babel}
\usepackage{amssymb,amsthm,amsmath,latexsym,amsfonts,
amscd,epsfig}
\usepackage{subfigure}

\usepackage{float}
\usepackage[section]{placeins}
\usepackage{pstricks}
\usepackage{pstcol}

\newcommand{\sect}[1]{\setcounter{equation}{0}\section{#1}}

\def\be{\begin{equation}}
\def\ee{\end{equation}}
\def\bea{\begin{eqnarray}}
\def\eea{\end{eqnarray}}
\def\R{{\mathbb R}}

\def\hbar{{\xi}}

\def\conm#1#2{\left [ {#1},{#2} \right ]}

\begin{document}

\begin{center}
{\LARGE{\bf{Bases in Lie and Quantum Algebras}}}
\end{center}

\bigskip\bigskip

\begin{center}
A. Ballesteros$^1$, E. Celeghini$^2$  and M.A. del Olmo$^3$
\end{center}

\begin{center}
$^1${\sl Departamento de F\'{\i}sica, Universidad de Burgos, \\
E-09006, Burgos, Spain.}\\
\medskip

$^2${\sl Departimento di Fisica, Universit\'a  di Firenze and
INFN--Sezione di
Firenze \\
I50019 Sesto Fiorentino,  Firenze, Italy}\\
\medskip

$^3${\sl Departamento de F\'{\i}sica Te\'orica, Universidad de
Valladolid, \\
E-47011, Valladolid, Spain.}\\
\medskip

{e-mail: angelb@ubu.es, celeghini@fi.infn.it, olmo@fta.uva.es}
\end{center}
\bigskip

\begin{abstract}

Applications of algebras in physics are related to the
connection of measurable observables to relevant elements of the
algebras, usually the generators. However, in the determination of
the generators in Lie algebras there is place for some arbitrary
conventions. The situation is  much more involved in the context of  quantum algebras, where
inside the quantum universal enveloping algebra, we have not
 enough primitive elements that allow for a privileged set  of generators  and all
basic sets are equivalent.   In this paper we discuss how  the Drinfeld double
structure underlying every simple Lie bialgebra characterizes
uniquely a particular basis without any freedom,    completing the Cartan
program on simple algebras. By means of a
perturbative construction, a  distinguished deformed basis (we call it the
analytical basis) is obtained for every quantum group  as
the analytical prolongation of the above defined Lie basis of the
corresponding Lie bialgebra.  It turns out that    the whole  construction
is unique,  so to each quantum universal enveloping algebra
is associated one and only one bialgebra.  In this way the problem of the classification of quantum algebras is
 moved  to the
classification of bialgebras.   In order to make this procedure more clear, we discuss in
detail the simple cases of $su(2)$ and $su_q(2)$.
\end{abstract}
\vskip 1cm

\vskip 0.4cm

\vfill
\eject


\sect{Introduction}

Algebras are, of course, the same   object  in mathematics and
physics. However, from the mathematical viewpoint the interest
is focused on the general description of the structure as abstract
as  possible, while for   physical applications the interest is   focused
on the basic elements  of the algebra  to which some precise physical meaning is attributed.

In this  ``physically" motivated work we summarize  a  research line devoted to
individuate, for both Lie and quantum algebras, these basic objects
--the generators-- in an intrinsic way   eliminating every
arbitrariness.  In this review we do not  describe all the technical details
of this approach and the interested reader is invited to look at the specific
papers quoted in the bibliography
\cite{ballesteros04}--\cite{CBO}.

In the Lie case the situation is  quite  simple and   almost solved
more than half  century ago. Inside the universal enveloping algebra   of a given Lie algebra
the  distinguished  subspace of the algebra is completely  determined by the
requirement to be closed under commutation relations or,
equivalently, to have a primitive coproduct. Also the decomposition
of the space of every simple Lie algebra in   the form  $n_- \oplus h \oplus
n_+$ (where $h$ is the abelian Cartan subalgebra and $n_\pm$ are two
nihilpotent equivalent subalgebras) is an intrinsic property  of
interest for both theory and applications. Then the Cartan analysis of
$n_+$ allows to individuate, always without any arbitrariness, the
simple roots and to draw with them the Dynkin diagram. However, some
freedom remains in the length of the  roots and different
conventions have been adopted  so far:   Weyl and
Chevalley canonical bases, etc.   (see, for instance, \cite{Cornwell}).
 In this context, a way to ``solve" the problem of non-simple roots  is just to
eliminate them from the game:  the Chevalley reformulation is
indeed made in terms of the Cartan subalgebra, the simple roots and
the Serre relations. However,   in Physics this is not the usual way to proceed
since, in general, all roots (either simple or not) are considered on the some footing and
are related to physical objects. In conclusion, few details    remained
 until now unsolved in the construction of the bases
within  the Cartan program on simple algebras.

As a first result, we discuss here how   these problems  can be solved
by  introducing  objects more general than algebras:  Lie bialgebras and
Drinfeld doubles. It is indeed well known that  every  simple Lie
algebra  has associated  a
Lie bialgebra which is  ``almost" a Drinfeld double \cite{CP}. This
``approximation" can be in fact removed  because
every simple algebra of rank $r$, ${\cal A}_r$, is related to the
Drinfeld double  ${\cal A}_r \oplus t_r$, where $t_r$ is an abelian
$r$-dimensional algebra \cite{ballesteros05,ballesteros06}.  It turns out that, up to  an   overall
multiplicative factor implicit in the algebra,  this
Drinfeld double fixes uniquely the normalization factor   for
all  the generators.

The situation is quite different  when we attempt to introduce a
quantum deformation. Indeed, a general investigation concerning the
uniqueness of this quantization process has not been given yet and
only some restrictive results for certain deformations of simple Lie
algebras have been obtained. As a consequence, although the
existence of the quantization is guaranteed for each Lie
bialgebra (see \cite{CP}), a   classification of
quantum groups in the spirit of Cartan cannot be found in the literature
\cite{ballesteros04,ballesteros05}.

Moreover, the infinitesimal counterpart of a quantum group   is not
well defined. So, we are not dealing with a deformation of a Lie
algebra but with a quantum universal enveloping algebra $U_q(g)$,
which is a Hopf algebra deformation of the universal enveloping
algebra $U(g)$ of the Lie algebra $g$, i.e. a deformation of an
infinite-dimensional object. A  possible basis of $U(g)$  is   the
Poincar\'e-Birkhoff-Witt (PBW) basis on $g$ (i.e. the set of all
ordered monomials of powers of the generators of $g$).  Of course,
there are infinitely many basic sets different from $g$ (all related
among them by invertible transformations) that originate suitable
PBW bases. However $g$ is a special subspace of $U(g)$ and, for this
reason the only one considered. But in quantum algebras it is not
easy to find a prescription that permits  the choice of a privileged
bases. Indeed, contrarily to the non deformed case, where inside all
the sets of basic elements the vector space of generators is
univocally defined, in $U_q(g)$ there is an infinitude of basic sets
coexisting on the same footing. This implies that quantum
deformations lead us locally (therefore, geometrically) to
structures which are quite different from Lie algebras. In
particular, when we consider a deformation of an $n$-dimensional Lie
algebra, only the infinite dimensional $U_q(g)$ makes sense in
despite that its PBW basis is constructed in terms of a basic set of
$n$ elements. This problem of the basis underlies many difficulties
encountered when a precise physical/geometrical meaning has to be
assigned to $U_q(g)$  as, for instance, in the context of quantum
deformations of space-time symmetries. In that case, it is well
known that the models so obtained depend on the choice of different
bases \cite{majid}, where different possibilities are related
through nonlinear transformations.

The second result   that we present here is the solution of the
above mentioned problem, providing a universal prescription for the
characterization and the explicit construction of a $n$-dimensional
basis, that can be properly called, the quantum analogue of the Lie
basis. To begin with,  we  analyze the role and properties of the
Lie algebra generators  within $U(g)$. Among the infinite possible
PBW bases, all related by nonlinear invertible transformations, the
generators determine the only one that is closed under commutation
rules and whose tensor product representations are constructed
additively. The latter property can be stated in Hopf-algebraic
terms as the Friedrichs theorem~\cite{jacobson}, which  states that
the only primitive elements in $U(g)$ ({\it i.e.} the elements such
that $\Delta(X)=\Delta_{(0)}(X):= 1\otimes X + X\otimes 1$) are just
the generators of $g$ as a Lie algebra. In this way, the generators
of $g$ become distinguished elements of $U(g)$. However, the
situation changes drastically in $U_q(g)$, where the law for the
construction of tensor product representations (coproduct) uses
nonlinear functions and no primitive bases exist. We show that,
among the infinitely many possible bases, there is just one where
the coproducts are ``as primitive as possible'', since all its
inessential terms, related to nonlinear transformations in the
enveloping algebra, have been removed. Thus, the only changes from
the Lie primitive coproducts are those imposed by the consistency
between the bialgebra cocommutator $\delta$ and the Hopf algebra
postulates. This basis (that has, of course, the same dimension than
the corresponding Lie algebra) is proposed as the true deformation
of the Lie algebra and  called ``quantum algebra",   to be
distinguished from the  (much larger) quantum universal enveloping
algebra.

To develop our approach we recall that Lie
group theory is based on analyticity with respect to group
parameters. In the same way, analyticity in the deformation
parameter(s) will give us the keystone for the identification of the
proper quantum algebra, that will be defined as the $n$-dimensional
vector space $(g_q,\Delta) \subset U_q(g)$ obtained as
analytical prolongation of the Lie bialgebra $(g,\Delta_{(0)})$.
Note that analyticity in the deformation has  already played a useful
role in quantum algebras, for instance in their contractions
\cite{LED}.

In this analytical prolongation, the cocommutator map of the Lie
bialgebra plays a fundamental role since it describes the first
order deformation and can be considered as the derivative at the
origin of the quantum coproduct. Thus, for a given Lie-Hopf algebra
$(g,\Delta_{(0)})$, there are as many different quantizations as
inequivalent cocommutators $\delta$. The Hopf algebra postulates
together with the assumption that all the results are analytical
functions of both the deformation parameter(s), introduced by the
bialgebra, and of generators determine uniquely the ``analytical"
generators.

As a third result we  will discuss how the Friedrichs theorem can be
extended to quantum algebras. Like in the Lie case, where to a
$U(g)$ corresponds one and only one Lie algebra $(g,\Delta_{(0)})$,
to a $U_q(g)$ corresponds one and only one quantum algebra
$(g_q,\Delta)$. Since we prove that $(g_q,\Delta)$ is in one-to-one
correspondence with $(g,\delta)$, then  to each bialgebra
corresponds one and only one quantization: the classification of
quantum groups is reduced to the classification of Lie bialgebras.

Summarizing, we close the commutative diagram by adding the vertex
$(g_q,\Delta)$
\[ \begin{CD}
(g,\Delta_{(0)}) @<{\rm \hskip0.42cm Friedrichs\; theorem\hskip0.42cm}<< U(g)  \\
@VV{q}V @ VV{q}V  \\ (g_q,\Delta) @<{\rm Generalized\; Friedrichs\;
th.}<< U_q(g)
\end{CD}\]
The vertical lines represent the quantization procedure, and the
horizontal ones are related to the definition of the basic set of
generators of the universal enveloping algebra  and of its quantum
analogue.

The paper is organized as follows. Section~\ref{drinfelddouble} is
devoted to Lie bialgebras and Drinfeld doubles, and  after a short
description of these objects we  explicitly discuss the simplest
case of $su(2)$,  where the length of the root vector turns out to
be fixed  by the self-dual Drinfeld double structure. The
generalization to any simple Lie algebra is briefly discussed. In
Section \ref{analyticaldeformation} we describe   the analytical
approach to quantization,  which allows to  relate,  in a
perturbative way,   the Lie bialgebra $(g,\delta)$  and  the quantum
algebra $(g_q,\Delta)$. In this way, we show that quantization not
only relates $U_q(g)$ to $U(g)$ but also quite simpler
$n$-dimensional objects like $(g_q,\Delta$) to $(g,\delta)$. Again
the discussion will be illustrated through the standard deformation
of $su(2)$. Starting from the standard bialgebra of $su(2)$ we
derive a well defined form of the quantum algebra $su_q(2)$.
Section~\ref{perturbativefriedrichstheorem} is devoted to revisit
the first horizontal line of the diagram. Its generalization, that
corresponds to the lower line, is    discussed in
Section~\ref{unicityofdeformation}. Finally,  some conclusions close
the paper.


\sect{Lie bialgebras and Drinfeld doubles}
\label{drinfelddouble}

A Lie bialgebra $(g,\eta)$ is a Lie algebra $g := \{ Z_p\}$
endowed with  a skew-symmetric linear map (cocommutator) $\eta :
g\to g\otimes g$ such that
\[
\eta([Z_p,Z_q])=({\rm ad}_{Z_p}\otimes 1 +1\otimes {\rm ad}_{Z_p} )
\eta (Z_q) - ({\rm ad}_{Z_q}\otimes 1 +1\otimes {\rm ad}_{Z_q} )
\eta (Z_p), \qquad \forall Z_p,Z_q \in g.
\]
Thus,  besides the commutation relations
 \be [Z_p,Z_q]=
f^r_{p,q}\,Z_r \, , \ee  we have another set of relations
 \be
\label{deltaZ} \eta(Z_p)= c_p^{q,r}\,Z_q\otimes Z_r . \ee The two
sets of structure constants $f^r_{p,q}$ and $c_p^{r,q}$ look quite
  symmetric and, indeed, on the dual space $g^*:= \{ z^p\}$ they
play an inverted role
 \be\label{deltaz}
 [z^p,z^q]=
c_r^{p,q}\,z^r, \qquad \eta(z^p)= f^p_{q,r}\,z^q\otimes z^r \,.\ee
The compatibility relations  between the structure constants
\be\label{compatibility}
c^{p,q}_r f^r_{s,t} = c^{p,r}_s f^q_{r,t}+c^{r,q}_s f^p_{r,t}
+c^{p,r}_t f^q_{s,r} +c^{r,q}_t f^p_{s,r}   \ee
allow to combine
$g$ and $g^*$ in an unique algebra called a  Drinfeld double 
\cite{CP}. This Drinfeld double is, indeed, defined as a Lie algebra
$\bar g$ such that,
 $\bar g = g \oplus g^*$ as vector space, where $g$ and
$g^*$ are disjoint sub-bialgebras, and the crossed commutation rules
between $g$ and $g^*$ are defined in terms of the structure tensors
$f^p_{q,r}$ and $c_p^{q,r}$ by eq. (\ref{zz}).
The Jacobi identities  for the algebra $\bar g$  are just   the
compatibility relations (\ref{compatibility}).

The crossed commutation relations \be [z^p,Z_q]= f^p_{q,r}z^r-
c^{p,r}_q Z_r . \label{zz} \ee are related to the invariance of the
pairing \be \label{pairingz} \langle Z_p,Z_q\rangle=0 ,\qquad
\langle Z_p,z^q\rangle=\delta_p^q , \qquad \langle z^p,z^q\rangle=0
. \ee

 Moreover, the   double  Lie algebra $\bar g$ has a (quasitriangular) double Lie
bialgebra structure $(\bar g,\delta)$, where its cocommutators are
determined by its sub-bialgebras ($g,\eta$) and ($g^*,\eta$) as
\be\label{tt} \delta(Z_p)= -\eta(Z_p)= -c^{q,r}_p Z_q\otimes Z_r
,\qquad \delta(z^p)=\eta(z^p)=f^{p}_{q,r} z^q\otimes z^r . \ee

In the particular case that $c_p^{q,r}=-f^p_{q,r}$ \cite{gomez} the
Drinfeld double gets an additional symmetry, and we call it a
self-dual Drinfeld double. It is worthy noticing that the positive
and negative Borel subalgebras $b_\pm$ of any simple Lie algebra $g$
have this last property, but, as they have the Cartan subalgebra in
common, they cannot be identified {\bf  as} $g$ and $g^*$. Hence,
the classical Lie algebras  are ``almost" Drinfeld doubles
\cite{DrICM,CP}.

 In fact, this  problem has been solved with a slight extension of $\bar g$
considering the direct sum of $\bar g$ with an abelian algebra
$t_n$, where $n$ is the rank of $\bar g $ \cite{ballesteros05}.
The most elementary example of self-dual Drinfeld double structure
is that of the algebra $su(2)\oplus t_1$. We write $su(2)$ as
\be\label{su2conmutators}
 [J_3,J_\pm]=\pm J_\pm ,\qquad[J_+,J_-]=
\alpha J_3 , \ee where the most common physical convention is
$\alpha = 2$ , while   in mathematical contexts different
conventions and a different overall factor are used to write
$su(2):= \{h,f,g\}$   relations as $[h,f]=2f,\;[h,g]=-2g,\;[f,g]=h$.
In disagreement with both of them, from $su(2)\oplus t_1$, we  will
obtain $\alpha = 1$.

The two $su(2)$ Borel subalgebras $ s^\pm :=\{ J_3, J_\pm\}$ are not
disjoint (as it should be required by (\ref{pairingz})) but, if we
introduce a new abelian algebra $t_1 \equiv \R$ with generator $I$,
that commutes with $su(2)$ (i.e. $[I, \cdot ]=0$), we can define two
disjoint algebras $b_+ := \{Z_1, Z_2\}$ and $b_- :=\{ z^1,z^2 \}$
that allow us to define a true self-dual Drinfeld double. Let us
indeed define
\[
Z_1= \frac{1}{\sqrt{2}}(J_3 + i I), \qquad Z_2=J_+, \qquad \qquad
z^1=\frac{1}{\sqrt{2}}(J_3 - i I) ,\qquad z^2= J_- .
\]
The commutation relations imposed by $su(2)\oplus t_1$ inside $b_+$
and $b_-$ are, respectively
\[
[ Z_1, Z_2] =\frac{1}{\sqrt{2}}Z_2,\qquad [ z^1, z^2]
=-\frac{1}{\sqrt{2}}z^2 ,
\]
such that, the structure  constants are
\[
f^{2}_{1,2}=-f^{2}_{2,1}=
-c^{1,2}_{2}=c^{2,1}_{2} =\frac{1}{\sqrt{2}}  .
\label{tensores}
\]
Compatibility conditions (\ref{compatibility}) are easily checked.
Hence, because $c^{i,j}_{k}= -f^{k}_{i,j}$, we have a self-dual
Drinfel double. From (\ref{zz}) we get the crossed commutation
relations  as
\[
[ Z_1,z^1] =0,  \qquad [ Z_1,z^2] =-\frac{1}{\sqrt{2}} z^2,  \qquad
[ Z_2,z^1] =-\frac{1}{\sqrt{2}} Z_2,  \qquad [ Z_2,z^2] =
\frac{1}{\sqrt{2}} (Z_1+z^1).
\]
The cocommutators are
\[
\delta(Z_1)=0 , \qquad \delta(Z_2)=\frac{1}{\sqrt{2}}Z_2\wedge Z_1
,\qquad \delta(z^1)=0 , \qquad \delta(z^2)=\frac{1}{\sqrt{2}}z^2
\wedge z^1 .
\]

Now, we can come back to our original basis. The commutation
relations inside $b_\pm$ have been, of course, assumed but the value
of $\alpha$ is determined by the Drinfeld double to be exactly
$\alpha=1$. This is equivalent to impose  the following relation on
the $so(3)$ generators
\[
J_\pm \equiv \frac{1}{\sqrt{2}} ( J_1 \pm i J_2 )
\]
as, from a different point of view, is required by the Killing form.

 In this way  the full double Lie bialgebra is obtained. Since we are interested in the
subalgebra $su(2)$, we consider the trivial representation of $t_1$
(where $I=0$) \be\label{delta_su2} \delta (I)=0 ,\qquad \delta
(J_3)=0 ,\qquad \delta (J_+)=\frac{1}{2}J_+\wedge J_3 ,\qquad \delta
(J_-)=\frac{1}{2}J_-\wedge J_3 . \ee

The generalization to larger algebras is straightforward. We display
here the results for the $A_n$ series ~\cite{ballesteros06}. Let us
consider  the Lie algebra $gl(n+1)=A_n\oplus h:= \{H_i,\; F_{ij}\}$,
where $h$ is the Lie algebra generated by $\sum H_i$, we define:
\be\begin{array}{l}\label{commutatorsA}
[H_i,H_j]=0, \\[0.3cm]
[H_i,F_{jk}]=(\delta_{ij} - \delta_{ik})F_{jk}, \\[0.3cm]
[F_{ij},F_{kl}]=(\delta_{jk} F_{il}- \delta_{il} F_{kj}) +
\delta_{jk} \delta_{il}(H_i - H_j) .
\end{array}\ee
As before, we introduce an abelian algebra $t_{n+1}:= \{I_{i}\}$
with $n+1$ central generators $I_i$ and we construct the algebra
$gl(n+1)\oplus t_{n+1}$ that we show to be a Drinfeld double. With a
change of basis similar to the $su(2)$ case, the new generators are
\be\label{change} Z_i:=\frac{1}{\sqrt{2}}(H_i + {\bf i} I_i), \quad
Z_{ij}:=F_{ij},\qquad z^i:=\frac{1}{\sqrt{2}}(H_i- {\bf i} I_i),
\quad z^{ij}:=F_{ji} \quad(i<j). \ee We obtain thus two equivalent
soluble Lie algebras  $b_\pm$ each of dimension $(n+1)(n+2)/2$
\[
b_+:=  \{ Z_i,Z_{ij}\} ,\qquad b_-:=   \{ z^i,z^{ij}\},\qquad
i,j=1,\dots,n+1,\quad i<j ,
\]
such that $b_+ + b_- = gl(n+1)\oplus t_{n+1}$ as vector spaces.

The  commutation rules imposed by (\ref{commutatorsA}) inside $b_+$
and $b_-$ are \be\begin{array}{llll}\label{snmas} & [Z_i,Z_j]=0,
\quad &[Z_i,Z_{jk}]=\frac{1}{\sqrt{2}}(\delta_{ij} - \delta_{ik})\,
Z_{jk},\quad
& [Z_{ij},Z_{kl}]=(\delta_{jk} Z_{il}- \delta_{il} Z_{kj}) ,\\[0.3cm]
& [z^i,z^j]=0, \quad & [z^i,z^{jk}]=- \frac{1}{\sqrt{2}}(\delta_{ij}
- \delta_{ik})\, z^{jk},\quad & [z^{ij},z^{kl}]=- (\delta_{jk}
z^{il}- \delta_{il} z^{kj}).
\end{array}\ee

The two algebras $b_+$ and $ b_-$ can be paired by
\be\label{pairingA} \langle z^i,Z_j\rangle= \delta^i_j ,\qquad
\langle z^{ij},F_{kl}\rangle= \delta^i_k\delta^j_l , \ee defining a
bilinear form such that both $b_\pm$ are isotropic.

The compatibility relations (\ref{compatibility}) can be easily
checked and we are now able to came back from
 (\ref{snmas}) and (\ref{pairingA}) to all the formulas (\ref{commutatorsA}) but
this is strictly related to our definition of the root vectors. A
change of scale like $F_{ij} \rightarrow k F_{ij}$ will introduce
into (\ref{commutatorsA}) changes that are is incompatible with the
Drinfeld double structure.

This procedure can be easily extended to any arbitrary semi-simple
Lie algebra $g_n$ by considering $g_n \oplus t_n$, which can be
equipped with a  self-dual Drinfeld double structure. In
Ref.~\cite{ballesteros06,ballesteros07}   the explicit form for the
classical series $A_n$, $B_n$, $C_n$ and $D_n$ can be found.

In conclusion: there is a Drinfeld double behind every simple Lie
algebra but only if it is written in the correct basis.

The Cartan program on semi-simple Lie algebras has been thus
concluded by fixing, without any freedom, the length of the vectors
associated to all the roots. This is the main result coming from the
Drinfeld double perspective at a ``classical" level.

 It is worthy to note  that the Drinfeld
double fixes everything not only at the  Lie algebraic level but
also at the  Lie bialgebra one: the canonical Lie bialgebra
structure for $gl(n+1)\oplus t_{n+1}$ is determined by the
cocommutator $\delta$  and reads
\be\begin{array}{ll}\label{cocommutadorA}
\delta (I_i)=0, \\[0.3cm]
\delta (H_i)=0,\\[0.3cm]
\delta (F_{ij})=-\frac 12 F_{ij} \wedge (H_i-H_j)- \frac {{\bf
i}}{2} F_{ij}
\wedge (I_i-I_j) + \sum_{k=i+1}^{j-1}{F_{ik} \wedge F_{kj}}, \qquad & i<j ,\\[0.3cm]
\delta (F_{ij})=\frac 12 F_{ij} \wedge (H_i-H_j)- \frac {{\bf i}}{2}
F_{ij} \wedge (I_i-I_j) - \sum_{k=j+1}^{i-1}{F_{ik} \wedge
F_{kj}},\qquad & i>j.
\end{array}\ee

We have to observe that the chain of Drinfeld doubles $g_n\oplus
t_n\subset g_{n+1}\oplus t_{n+1}$ is preserved at the level of Lie
bialgebras. However, although  $g_n$ is a subalgebra of $g_n\oplus
t_n$, the cocommutator $\delta (g_n)$ does not define a
sub-bialgebra since it depends on the extra $t_n$ sector.

\sect{Analytical deformation}\label{analyticaldeformation}

 The problem of the generalization of the previous results to quantum algebras is
quite more difficult.  Strictly speaking, quantum algebras are not deformed Lie
algebras in the sense that  the initial object   to be deformed is $U(g_n)$, the
(infinite-dimensional) universal enveloping
algebra built on the Lie algebra $g_n$, and the final object,
with similar properties, is the quantum universal enveloping algebra
$U_q(g_n)$ (also  infinite-dimensional).

The universal enveloping algebra $U(g_n)$  has, as a possible basis, the
Poincar\'e-Birkhoff-Witt (PBW) basis,  i.e. the set of all ordered
monomials built up on the powers of the generators of $g_n$.
Moreover, the PBW basis is only a particular one among the possible bases   all   related between
them by invertible transformations. The crucial  point is that, when
considering $U_q(g_n)$, the  algebra subset looses all its
privileges in a structure where linearity does not make sense.
On the contrary, when considering  $U(g_n)$, we have the chance that
the Lie subset is the only subset  that has two relevant
properties: first, it is the only $n$-dimensional set closed under
commutation and, second, only the elements of the Lie algebra are
primitive. This last property is --at least in physical terms-- the
most relevant as it implies the additivity of physical observables
related to the algebra.

Indeed,  $U(g_n)$ and $U_q(g_n)$ are equipped with  Hopf algebra
structures. The main point is the existence of an isomorphism,
called  coproduct,
\[
\Delta : U(g_n)\to U(g_n)\otimes U(g_n), \qquad
\Delta _q: U_q(g_n)\to U_q(g_n)\otimes U_q(g_n).
\]
Those elements of $U(g_n)$ such  that $\Delta (X)=1\otimes
X+X\otimes 1$ are said to be primitive and the Friedrichs theorem
establishes that  only the   elements of $g_n$ are
primitive. Moreover, as the Lie basis elements are cocommutative
(i.e. $\Delta (X)=\sigma\circ\Delta (X)$ with $\sigma (A\otimes
B)=B\otimes A$), this property is also valid for all $U(g_n)$.

Unfortunately,  deformation breaks both  properties,  primitivity and cocommutativity.
So, we looses  the canonical procedure to
individuate  among all the possible basic sets  a privileged one. At this point
 we are compelled to considerer at the same time the
whole $U_q(g_n)$  where  all the bases are equivalent. Let us recall
that three basic sets are usually considered in the literature:
 Lie basis, crystal  basis \cite{Lu,Ka} and the most common  one,  which
 for  $su_q(2) := \{J_3,J_\pm\}$ \cite{CP}  is
\be\label{basissuq2a}
 [J_3,J_\pm]=\pm J_\pm ,\qquad
 [J_+,J_-]= \frac{\sinh (zJ_3)}{\sinh (z/2)} ;
 \ee
 \be\label{basissuq2b}
\Delta(J_3)=1\otimes J_3+J_3\otimes 1,\qquad
\Delta(J_\pm)=e^{\frac{z}{2}J_3}\otimes J_\pm+J_\pm\otimes
e^{-\frac{z}{2}J_3}.
 \ee

 The Lie basis of $U_q(g_n)$ describes a situation where both the basic set and its
coproduct are closed under commutation relations but the relation
between them is different (and quite more complex) from the Lie one.
For a physicist this means that we have the same symmetry of the Lie
situation but we loose the usual definition of composed system. It
could, perhaps, be the appropriate scheme to describe systems of
particles with an interaction among them. The crystal
 basis (often called canonical basis) is, in some sense, the opposite of
the Lie basis: while in the Lie basis the algebra is written for $z
\to 0$ and all deformation is transferred to the coproduct, in the
crystal basis the algebra is obtained for $z \to \infty$ while, as
in the Lie basis, the problem of consistency  is left to the
coproduct. Also in the crystal basis the representations look quite
simple and have been applied in statistical mechanics
\cite{Kang} and in genetics \cite{MS}.

However, it is not an accident that everybody knows   $U_q(su(2))$
in terms of the basic set  displayed in (\ref{basissuq2a}),
(\ref{basissuq2b}),  since this basis is, up to few details,  the analytical
prolongation to the $su(2)$ generators defined in
Sect.~\ref{drinfelddouble}. More information can be found in Ref.
\cite{CBO}, where the discussion is more general and technical.

Analyticity in Lie theory is analyticity in the parameter (dual
space of the generators), but here we have to extend the analyticity
to the deformation parameters $z$ in such a way that in a well precise limit (usually $z\to
0$)  the quantum algebra  reduces to the corresponding Lie algebra (without loss of generality,
we shall restrict here to the one-parameter case).   We also recall that analyticity in quantum
algebras is not new since, for instance, contraction theory relies upon this concept \cite{LED}.

The point is that from the theory of analytical functions we know
that there exist  infinitely many analytical functions that have the same
value in a point and the same derivative in the same point.
Consistently, from a Lie bialgebra we can obtain, by analytic
prolongation,   infinitely many basic sets all of them describing the same
$U_q(g)$ and expressed by means of  infinitely many free parameters,
independent from $z$. This is coherent with the fact that $U_q(g)$
(as well as $U(g)$) can be described by  infinitely many basic sets.
But, as discussed in the following, all these free parameters remain
present in the limit $z\to 0$: each one of the basic sets of
$U_q(g)$ is related by analyticity to a basic set of $U(g)$ and,  as
$U(g)$ has a privileged basic set in the Lie generators, the basic
set related by analyticity to Lie generators is a privileged basic
set for $U_q(g)$.

For instance $U(su(2))$ can be defined, as usual, from the $su(2)$
generators, but --equivalently-- starting from another basic set
like, for instance, \be\label{cocomm} Y_\pm = J_\pm + J_3^2 J_\pm ,
\qquad Y_3 = J_3 + J_3^2  .\ee As these relations are (formally)
invertible $U(su(2))$ --as well as its deformation-- does not
change, but $\{Y_\pm, Y_3\}$ is clearly not the best choice for
applications. This is exactly the freedom we remove in the
definition of the analytical basis  by imposing that the deformation
starts from the Lie generators and that only the contributions
necessary to save the consistency between the bialgebra ($g,\delta$)
and the Hopf postulates are introduced. All contribution related to
nonlinear transformations in $U(g)$, like (\ref{cocomm}),  are in
this way eliminated.

As it is well-known $\Delta$ is related to $\delta$ by \be\label{dD}
\delta =\lim_{z\to 0}\frac {\Delta -\sigma\circ \Delta}{2 z} :
 \ee
the cocommutator   $\delta$ can be seen as the derivative at the
origin of the quantization and $U_q(g)$ is sometimes called a
``quantization of $U(g)$ in the direction of $\delta$''.

The analytical deformation is introduced in three steps:
\begin{enumerate}
\item We find order-by-order the changes imposed by
$\delta$ in $\Delta_{(0)}$  to be consistent with coassociativity
(\ref{recursionk}) and we determine in this way the full coproduct $\Delta$.
\item
 By using analyticity and the homomorphism property of $\Delta$
 (\ref{homomorphismproperty}), we obtain the
 commutation rules for $g_q$ starting from the known ones of $g$. Thus, the
 $n$-dimensional $(g_q,\Delta)$ is constructed.
\item
From $g_q$ a PBW basis in $U_q(g)$ is built.
\end{enumerate}

In this way we construct an unique connection $(g,\delta)
\rightarrow (g_q,\Delta) \rightarrow U_q(g)$. Since $(g,\delta)$ is
the limit of $(g_q,\Delta)$ and, as it will be show in
Sect.~\ref{unicityofdeformation}, every $U_q(g)$ admits one and only
one basis $g_q$, the arrows can be inverted and a one-to-one
correspondence is found between $(g,\delta)$ and $U_q(g)$. So,
equivalences between $q$-deformations of  $U_(g)$ imply equivalences
among their associated  bialgebras and the classification of the
$U_q(g)$ is carried to the quite simpler classification of Lie
bialgebras.

Analyticity means that the commutation relations of any basic set
$\{Y_j\}$ $(j=1,2,\dots, n)$ of $U_q(g)$ (as well as of $U(g)$) are
analytical functions of the $Y_j$ and that the quantum coproduct
$\Delta$ of the $Y_j$ can be written as a formal series
\be\label{coproductseries} \Delta(Y_i)=
\sum_{k=0}^{\infty}\Delta_{(k)}(Y_i)= \Delta_{(0)}(Y_i) +
\Delta_{(1)}(Y_i)  + \dots \ee where $\Delta_{(k)}(Y_i)$ is a
homogeneous polynomial of degree $k+1$ in \,$1\otimes Y_j$\, and
\,$Y_j \otimes 1$.

There are two properties of Hopf algebras   which are relevant in this analytical
approach: the coassociativity condition \be
\label{coass} (\Delta\otimes 1-1\otimes\Delta )\circ\Delta (Y_i)=0,
\ee that in a perturbative form will be rewritten as
\be\label{recursionk}
 \sum_{j=0}^{k}\left( \Delta_{(j)}\otimes 1 - 1
\otimes \Delta_{(j)} \right) \circ \Delta_{(k-j)} (Y_i)=0
,\qquad\qquad \forall k ;\ee and the homomorphism property
\be\label{homomorphismproperty1} \Delta(\conm{Y_i}{Y_j})=
\conm{\Delta(Y_i)}{\Delta(Y_j)} ,\ee which can be rewritten as
\be\label{homomorphismproperty} \Delta_{(k)}(\conm{Y_i}{Y_j})=
\sum_{l=0}^k \conm{\Delta_{(l)}(Y_i)}{\Delta_{(k-l)}(Y_j)}
,\qquad\qquad \forall k \,. \ee

To enlighten the construction we discuss explicitly the
quantization of the  Lie bialgebra $(su(2),\delta)$ described by
(\ref{su2conmutators}) with $\alpha =1$ and
(\ref{delta_su2}).
More details can be found in \cite{CBO}. The results of this Section
show that analyticity chooses the coproduct given by formula
(\ref{basissuq2b})  among all possible ones \cite{ballesteros04} and
the usual expressions (\ref{basissuq2a}) must be   replaced by other
slightly different  ones.

As stated before,  let us start by getting the quantum coproduct. The
zero-approximation of the quantum coproduct is $\Delta_{(0)}$, while
$\Delta_{(1)}$ is essentially $\delta$ and the
$\Delta_{(k)}$ are obtained adding order by order the contributions
imposed by  consistency with eq. (\ref{recursionk}).

The case of $J_3$ is simple: $\Delta_{(0)}(J_3) = J_3 \otimes 1 +
1\otimes J_3$ and $\delta(J_3)=0$, that implies that the
anti-cocommutative part of $\Delta_{(1)}(J_3)$ is zero. Thus
$\Delta_{(1)}(J_3)$ must be cocommutative. But cocommutative
contributions are related to non linear bases like the one of eq.
(\ref{cocomm}) in the Lie algebra and can be removed by coming back,
with an appropriate change of basis, to the Lie generators. So (see
eq.(\ref{0to1}) and Ref.~\cite{CBO} for an exhaustive discussion) we write
$\Delta_{(1)}(J_3)=0$. Then, we have
 \[
 \Delta(J_3)=
\Delta_{(0)}(J_3) + {\cal O}_{(2)}(J_3), \] where we have defined
${\cal O}_{(m)}(Y_i)$ as a polynomial or degree $> m$ in $1\otimes
Y_j$ and $Y_j\otimes 1$. From eq.(\ref{coproductseries}),
 we can write
\[\Delta(J_3)=
\Delta_{(0)}(J_3) + \Delta_{(2)}(J_3)+ {\cal O}_{(3)}(J_3).
 \]
Again, as eq. (\ref{recursionk}) for $k=2$ is  consistent with
$\Delta_{(2)}(J_3)=0$, we perform a (cubic) change of
basis and we write
 \[\Delta(J_3)=
\Delta_{(0)}(J_3) + \Delta_{(3)}(J_3)+ {\cal O}_{(4)}(J_3) .
\]
As this procedure can be iterated for all $k$, the analytical prescription for $J_3$ imposes
$\Delta_{(k)}(J_3)=0,\;\; \forall k>0$. Hence
\be\label{DE}\Delta(J_3) =\Delta_{(0)}(J_3) =J_3 \otimes 1 +
1\otimes J_3 , \ee i.e., the analytical procedure associates a
primitive coproduct
 to any null $\delta$. This could
be considered trivial, as $\{J_3\}$ closes an $u(1)$ algebra, but  this is not the case since
it has been obtained from a precise prescription and
it is not, like in \cite{CP}, an arbitrary choice inside $U(g)$.

For $\Delta(J_+)$ we have
\[
\Delta(J_+)= \Delta_{(0)}(J_+)+\Delta_{(1)}(J_+) + {\cal O}_{(2)}(J_+) .
\]
From (\ref{dD}) and the coassociativity condition (\ref{recursionk})
for $k=1$, we have that  $\Delta_{(1)}(J_+) = \frac z2
\,\delta(J_+)$ as, like in the $J_3$ case, the possible
cocommutative contribution  --independent from $\delta$--  can be
put to zero by a change of basis. We have thus
\[
\Delta(J_+)= \Delta_{(0)}(J_+)+\frac z2 \,\delta(J_+) + \Delta_{(2)}(J_+)+{\cal O}_{(3)}(J_+) .
\]
The coassociativity condition (\ref{recursionk}) for $k=2$ solved in
the unknown $\Delta_{(2)}(J_+)$ gives (again disregarding arbitrary
contributions that can be removed by a nonlinear change of basis)
 \be\label{D2X}
 \Delta_{(2)}(J_+)= \frac{z^2}{8}
 (J_3^{\; 2} \otimes J_+ + J_+ \otimes J_3^{\; 2}) .
 \ee
Considering next order, we can now write
\[
\Delta(J_+)= \Delta_{(0)}(J_+)+\frac z2 \,\delta(J_+) +
\Delta_{(2)}(J_+)+ \Delta_{(3)}(J_+)+{\cal O}_{(4)}(J_+) ,
\]
where $\Delta_{(2)}(J_+)$ is given by eq. (\ref{D2X}) and
$\Delta_{(3)}$ is the new unknown. Solving the coassociativity
condition (\ref{recursionk}) for $k=3$ we have
 \be\label{D3X}
 \Delta_{(3)}(J_+)= \frac{z^3}{48} (J_3^{\;3} \otimes J_+ - J_+ \otimes J_3^{\; 3})
 \ee
and the general formula obtained by iteration, always putting to
zero the cocommutative contributions related to nonlinearity in
$U(su(2))$ and independent from $\delta$, is

\[
 \Delta_{(k)}(J_+)= \frac{1}{k!}\left( \frac{z}{2}\right) ^k (J_3^{\; k} \otimes J_+
 +(-1)^k  J_+ \otimes J_3^{\; k}) ,\qquad
 \forall k .
 \]

\noindent The $\Delta_{(k)}$ are easily summed to
\[
 \Delta(J_+)= e^{\frac z2 J_3} \otimes J_+ + J_+ \otimes e^{-\frac z2 J_3} .
 \]

The approach is exactly the same for $J_-$ and gives a similar result.
Thus, we obtain the analytical quantum coproduct associated to
$(su(2),\delta)$
 \be
\label{su2cuantizacion1}\begin{array}{l}
 \Delta (J_3) = J_3 \otimes 1 + 1 \otimes J_3 ,\\[0.3cm]
 \Delta (J_+) = e^{\frac z2\,J_3} \otimes J_+
 + J_+ \otimes e^{-\frac z2\,J_3}, \\[0.3cm]
 \Delta (J_-) = e^{\frac z2\,J_3} \otimes J_-
 + J_- \otimes e^{-\frac z2\,J_3}  .
\end{array}\ee
 that it is exactly the coproduct displayed in expression
 (\ref{basissuq2b}).

 Now we have simply to  start from the
commutators (\ref{su2conmutators}) with $\alpha=1$  and
   to impose order by order the
homomorphism condition for the deformed commutation rules. We thus
find
 \be \label{su2cuantizacion2} [J_3,J_\pm]= \pm J_\pm ,\qquad
 [J_+,J_-]=\frac{1}{z}\,\sinh(z\,J_3) .
 \ee
 Note that $\Delta_{(0)}\left([J_+,J_-]\right) =
\Delta_{(0)}(J_3)$  forbids the introduction of $z$-dependent
normalizations, like $\sinh (z/2)$ of expressions
(\ref{basissuq2a}), and a factor 2 has been removed at the classical
level from the Drinfeld double structure. Expressions
(\ref{su2cuantizacion1}) and (\ref{su2cuantizacion2}) define
uniquely the standard analytical deformation of the Cartan basis of
$su(2)$ such that the $q$-generators $J_3,J_+,J_-$ could be called
the $q$-Cartan basis of $su_q(2)$.

By inspection, all the symmetries (for example, $\{J_3,J_+,J_-\}
\leftrightarrow \{J_3,-J_+,-J_-\}$) and the embedding conditions
(for instance, $su(2) \supset borel(J_3,J_+) \supset u(1):= \{
J_3\}$) of the bialgebra $(g,\delta)$ are automatically preserved in
the quantization $(g_q,\Delta)$. Like for the bialgebra only the
$u(1)$ generated by $J_3$ is a sub-quantum-algebra of $su_q(2)$.


\sect{Perturbative Friedrichs theorem}
\label{perturbativefriedrichstheorem}

We give now a constructive proof of the Friedrichs theorem,
building explicitly the primitive generators  $\{X_j\}$ in terms of
an arbitrary set $\{Y_j\}$ on which a PWB basis for the whole $U(g)$
can be built. The elements $Y_j$ are cocommutative but, in
principle, non primitive. The machinery consists in repeated changes
of bases that allow to obtain each time  a better approximation to
primitivity where --this is the essential point-- the problem is
reformulated at each step in terms of the preceding basis. One
infinite iteration of the procedure allows to find,  among the
infinitely many possible bases of $U(g)$, the Lie generators.

In more detail, let us consider $X_i \equiv \lim_{k\to \infty} X_i^k $ were
$\{X_i^k\}$ is a basic set that approximates the Lie-Hopf coproducts
up to order $k$. The   terms ${\cal O}_{(m)}(X_i^k)$  are, in this Section,
cocommutative since we are working
in $U(g)$. Now any original basic set $\{Y_i\}$ is a zero
approximation to $\{X_i\}$, i.e.  $X_i^0 := Y_i$. Indeed \be\label{aa}
 \Delta(X_i^0)
= \Delta_{(0)}(X_i^0) + {\cal O}_{(1)}(X_i^0) = X_i^0\otimes 1 +
1\otimes X_i^0 + {\cal O}_{(1)}(X_i^0). \ee
The explicit form consistent with formula (\ref{recursionk}) of ${\cal
O}_{(1)}(X_i^0)$ in expression (\ref{aa}) is
\[\label{O11}
  {\cal O}_{(1)}(X_i^0)
= \sum d_i^{jl}\; \left(X_j^0\otimes X_l^0 + X_l^0\otimes
X_j^0\right)
 \;+ \;{\cal O}_{(2)}(X_i^0),
\]
where $d_i^{lj}$ are arbitrary constants. From
expression (\ref{coproductseries}) we get
 \be\label{O1}
  {\cal O}_{(1)}(X_i^0)= \Delta_{(1)}(X_i^0)+
{\cal O}_{(2)}(X_i^0),\ee
but if we define (see \cite{CBO} for details) the next approximation
\be \label{0to1}
X_i^1 := X_i^0 -\;\sum d_i^{jl}\; \{X_j^0, X_l^0\}
\ee
we obtain  a coproduct for $X_i^1$ with   vanishing
first order contributions, i.e.
 \[
 \Delta(X_i^1) =
\Delta_{(0)}(X_i^1)+ {\cal O}_{(2)}(X_i^0). \]
Eq. (\ref{0to1}) allows to rewrite
${\cal O}_{(2)}(X_i^0)$   in terms of $X_i^1$ as
${\cal O}_{(2)}(X_i^1)$
 \[ \Delta(X_i^1) = \Delta_{(0)}(X_i^1) + {\cal O}_{(2)}(X_i^1). \]
Now, defining
 a new change of basis $ X_i^2 $  we
obtain, in terms of the same $X_i^2$, the coproduct of the second approximation $X_i^2$ to the
generators as
\[
\Delta(X_i^2) = \Delta_{(0)}(X_i^2) + {\cal O}_{(3)}(X_i^2) ,
\]
which is  free form both first and second order contributions.
The procedure can now be iterated and the
$\Delta_{m}(X_i^{m-1})$ contribution eliminated through a new change
of basis that affects the higher orders only. The residual term
becomes ${\cal O}_{(m+1)}(X_i^m)$ and we get the $m$-order approximation to
the Lie generators
\[
 \Delta(X_i^m) = \Delta_{(0)}(X_i^m) + {\cal O}_{(m+1)}(X_i^m).
\]

The true generators of the Lie algebra $g$ are (formally) recovered
in the limit \[ X_i := \lim_{m\to \infty} X_i^m \] and, in agreement
with the Friedrichs theorem, their coproduct is the primitive one \[
\lim_{m\to \infty} \Delta(X_i^m) =\lim_{m\to \infty}
\Delta_{(0)}(X_i^m) =\Delta_{(0)}(X_i) = \Delta(X_i) = X_i\otimes 1
+ 1\otimes X_i. \]
 This coproduct is an
algebra homomorphism with respect to the  Lie
commutation rules
 \[ \Delta[X_i, X_j] = [X_i,X_j]\otimes 1 +
1\otimes [X_i,X_j] ,\] and the $n$  generators of the Lie algebra are
univocally identified in a constructive manner within $U(g)$, pushing
away order by order the corrections to a primitive coproduct.

The central point of this perturbative
approach to Friedrichs theorem (as well as to its following
extension to $U_q(g)$) is that at each order all the relations can be
rewritten in terms of the corresponding approximations of the
generators.


\sect{Unicity of deformation}\label{unicityofdeformation}

This novel proof of Friedrich's Theorem has been described here
because the procedure that allows to individuate the generators of
the quantum algebra $g_q$, among the infinitely many possible bases
within $U_q(g)$, is exactly the same that allows to individuate the
generators of the Lie algebra $g$ among the infinitely many possible
bases of $U(g)$. The perturbative approach allows thus a
generalization of the Friedrichs theorem to quantum algebras.

Indeed, the construction introduced in
Sect.~\ref{perturbativefriedrichstheorem} works also when
$\delta\neq 0$. Hence,  it provides a prescription for the
construction of the almost primitive generators starting from an
arbitrary set of basic elements of any $U_q(g)$. In this way we
obtain $(g_q,\Delta)$ from $U_q(g)$ and we  close the diagram
displayed in the Introduction.

Let us sketch the proof of such result. As in the undeformed case,
let $\{Y_j\}$ be an arbitrary set of basic elements
 that determine  $U_q(g)$,  whose
classical limit is a Lie bialgebra with $\delta \neq 0$.
Eqs. (\ref{aa}) and (\ref{O1}) are still valid, but now we have to add to the cocommutative part of $\Delta_{(1)}$ the anti-cocommutative contribution given by
$\delta$
\[ \Delta(X_i^0)= \Delta_{(0)}(X_i^0)+ z \,\delta(X_i^0) + \sum
d_i^{jl}\; \left(X_j^0\otimes X_l^0 + X_l^0\otimes X_j^0\right) \; +
{\cal O}_{(2)}(X_i^0). \] Like in
Sect.~\ref{perturbativefriedrichstheorem} we define $X_i^1 := X_i^0
-\;\sum d_i^{jl}\; \{X_j^0, X_l^0\}$. As this change of variables
does affect the $\delta$ contribution only to higher orders, the
differences between $\delta(X_i^0)$ and $\delta(X_i^1)$ can be
included in ${\cal O}_{(2)}(X_i^0)$ (or, equivalently, ${\cal
O}_{(2)}(X_i^1)$). So,
\[ \Delta(X_i^1) = \Delta_{(0)}(X_i^1) + z\, \delta(X_i^1) +
{\cal O}_{(2)}(X_i^1). \]
From  eq. (\ref{coproductseries}),
 we introduce
$\Delta_{(2)}(X_i^{1})$
 \[
{\cal O}_{(2)}(X_i^{1})= \Delta_{(2)}(X_i^{1})+ {\cal
O}_{(3)}(X_i^{1}). \] Similarly to $\Delta_{(1)}$, the second order
$\Delta_{(2)}$ contains two contributions: the first one
--characteristic of the deformations-- is proportional to $z^2$, and
is constrained by the consistency between $\delta$ and the
coassociativity condition. The second one is only related to the
definition of the basis in $U(g)$ and can be removed by another
change of basis that does not modify the form of the $z$-dependent terms
since the modifications can be included in ${\cal O}_{(3)}(X_i^2)$.

This procedure can be iterated. Thus, once the problem is solved for
$\Delta_{(m-1)}$, a contribution to  $\Delta_{(m)}$ proportional to
$z^m$ (cocommutative for $m$ even and anti-cocommutative for $m$
odd) is found. Indeed the unessential $z$-independent terms are
removed, exactly as in the case $\delta = 0$, with a change of basis
that does not affect the form of the previous $z$-depending terms
because the introduced changes are always of orders higher of $z^m$
and thus pushed out in ${\cal O}_{(m+1)}$. For $m \to \infty$ the
same coproducts derived from $(g,\Delta_{(0)})$ in
Sect.~\ref{analyticaldeformation} are found. The deformed
commutation rules are then imposed by the homomorphism
(\ref{homomorphismproperty}). Hence, we have closed the lower row of
the diagram, finding that $U_q(g)$ has as one of its basic sets the
analytical basis $(g_q,\Delta)$ obtained as the analytical
continuation in $z$ of the Lie generators.


\sect{Concluding remarks}

 The main result of this paper is the definition and the operative prescription
for the complete construction of the basis in semi-simple algebras,
both Lie and quantum. In the Lie case our result end the Cartan program by fixing the
few normalizations  adjusted before by convention.
The situation is completely different in quantum
algebras where also the basic sub-space of $U_q(g)$, corresponding
to the Lie algebra, has to be found. We solve the problem by means
of analyticity: by starting from the cocommutator $\delta$ given by
the Lie bialgebra, a perturbative method based on cocommutativity
and followed by the appropriate sum of the series allows us to
determine the full expressions for the deformed coproducts. In this
way, each primitive coproduct is deformed to its quantum
counterpart. A perturbative application of the homomorphims property
of the coproduct map allows to start from the Lie form of the
commutators and to arrive to their deformations.
In this way the full analytical quantization --including non simple
roots-- is obtained.

As an important consequence,  this quantization procedure imposes
that the deformation of the commutators remain antisymmetric, i.e. in
the analytical basis there is no room for $q$-commutators,
objects without a precise symmetry. This is relevant for physical
applications since the commutators are essential elements in Quantum
Mechanics as well as in Poisson-Lie structures (whose
$q$-counterparts are not well-defined).

The realization of one and only one quantum analog of the generators
has also a relevant byproduct. As $(g,\delta) \leftrightarrows
(g_q,\Delta) \leftrightarrows U_q(g)$, a one-to-one correspondence
between $U_q(g)$ and $(g,\delta)$ is established. In other words,
each Lie bialgebra admits one (as it was known) and only one (as
stated here) quantization and the problem of classification of
semisimple quantum groups can be fully solved by restoring to the
classification of the underlying Lie bialgebras.

Finally let us note that the proposed method is constructive and it
could be implemented to build by computer the quantum algebra
deformation of every Lie bialgebra $(g,\delta)$.

\section*{Acknowledgments}
This work was partially supported  by the Ministerio de Educaci\'on
y Ciencia  of Spain (Projects FIS2005-03959 and MTM2007-67389), by
the Junta de Castilla y Le\'on (Project VA013C05) and by INFN-CICyT
(Italy-Spain).




\end{document}